# Option pricing in time-changed Lévy models with compound Poisson jumps

**Roman V. Ivanov**[a,*], **Katsunori Ano**[b]

[a]*Laboratory of Control under Incomplete Information,
Trapeznikov Institute of Control Sciences of RAS,
Profsoyuznaya 65, 117997 Moscow, Russian Federation*
[b]*Department of Mathematical Sciences,
Shibaura Institute of Technology,
Tokyo 135-8548, Japan*

roivanov@yahoo.com (R.V. Ivanov), kano2@mac.com (K. Ano)

**Abstract** The problem of European-style option pricing in time-changed Lévy models in the presence of compound Poisson jumps is considered. These jumps relate to sudden large drops in stock prices induced by political or economical hits. As the time-changed Lévy models, the variance-gamma and the normal-inverse Gaussian models are discussed. Exact formulas are given for the price of digital asset-or-nothing call option on extra asset in foreign currency. The prices of simpler options can be derived as corollaries of our results and examples are presented. Various types of dependencies between stock prices are mentioned.

## 1 Introduction

In recent years, more realistic models than the classic Brownian motion for the specification of financial markets were suggested and investigated. The generalized hyper-bolic distributions were introduced in Barndorff-Nielsen [3]. These distributions are infinitely divisible and hence generate a particular class of Lévy processes which can be represented as time-changed Brownian motions, see the monographs by Barndorff-Nielsen and Shiryaev [6] or Cont and Tankov [9] for details. Another important class is the generalized tempered stable distributions which were firstly introduced in Koponen [23] and then investigated in particular by Bianchi et al. [8] and Rosinski [37]. The generalized tempered stable processes are not always time-changed Brownian motions (see Küchler and Tappe [24] or Küchler and Tappe [25] on the bilateral gamma processes), although for example the variance-gamma process and CGMY process can be decomposed in this way. We refer on these facts to Madan et al. [33] and Madan and Yor [32], respectively.

---

[*]Corresponding author.

The variance-gamma process is the one of the most popular examples of the generalized tempered stable processes. The variance-gamma distribution was firstly proposed as a model for financial market data in Madan and Seneta [31] and Madan and Milne [30]. They discussed the symmetric case of the distribution. The properties of the variance-gamma process defined as the time-changed by gamma subordinator Brownian motion with drift were considered in Madan et al. [33]. Also, Madan et al. [33] gave the analytical expression for the European call option price in the variance-gamma model together with the definition of the process as the difference of two gamma ones. Further, a number of papers confirmed statistically the idea of using the variance-gamma process for the modeling financial indexes. Daal and Madan [10] and Finlay and Seneta [14] approved the variance-gamma model for the currency option pricing and the exchange rate modeling. Linders and Stassen [26], Moosbrucker [34] and Rathgeber et al. [36] simulated by the variance-gamma distribution the Dow Jones index returns. Mozumder et al. [35] considered the S&P500 index options in the variance-gamma model. Luciano and Schoutens [27] modeled the S&P500, the Nikkei225 and the Eurostoxx50 financial indexes by the variance-gamma process. Luciano et al. [29] and Wallmeier and Diethelm [43] confirmed the using of variance-gamma distribution for the modeling of the US and the Swiss stock markets, respectively.

The normal-inverse Gaussian distribution was introduced in Barndorff-Nielsen [3] to model some facts in geology as a member of the class of generalized hyperbolic distributions. Financial market data, including the Danish and the German ones, was specified then by the normal-inverse Gaussian process in Barndorff-Nielsen [4] and Rydberg [38]. Properties of the normal-inverse Gaussian process discussed as the time-changed by inverse-Gaussian subordinator Brownian motion were considered in Barndorff-Nielsen [5] and Shiryaev [41]. The normal-inverse Gaussian distribution in the context of risk modeling was discussed in Aas et al. [1] and Ivanov and Temnov [21]. Figueroa et al. [13] showed that the normal-inverse Gaussian distribution specifies well a high frequency data from the US equity markets. Teneng [42] proved that the normal-inverse Gaussian process fits to the dynamics of many various foreign exchange rates. Göncü et al. [15] confirmed that this distribution also relates to the statistics of emerging market stock indexes. The modeling of Bloomberg closing prices by the variance-gamma and the normal-inverse Gaussian distributions was discussed in Luciano and Semeraro [28].

If we discuss the problem of computing in Lévy models, the basic method is the Fourier transform one, see for details the review paper by Eberlein [11]. However, it puts some restrictions on the properties of the process or the type of the derivative payoffs. In particular, it can be shown that this method cannot be applied to the pricing of digital options in the volatile variance-gamma model or in the normal-inverse Gaussian one. The method of closed form solutions which had been introduced by Madan et al. [33] was proceeded then in the papers by Ivanov and Ano [20], Ivanov [18] and Ivanov [19] for the variance-gamma distribution and by Ivanov [17] and Ivanov and Temnov [21] for the normal-inverse Gaussian one. This paper continues the elements of the research by Madan et al. [33]. We discuss the problem of multi-asset digital option pricing in the variance-gamma model in the presence of extra downside compound Poisson jumps. These jumps reflect the influence of events which can evoke dramatic drops of assets on financial markets. The examples are the terror attack of 9/11, the Subprime mortgage crisis of 2007, the collapse of Lehman Brothers or the recent deep fall of oil prices. In Sections 3 and 4 the variance-gamma and the normal inverse-Gaussian models are considered, respectively. The obtained formulas give the option prices under different types of dependencies between the asset dynamics.

## 2 Setup and notations

We suggest that the risky asset log-returns $H_t^j = \log S_t^j$, $j = 1, 2, 3$, $t \leq T$, follow the sums of time-changed Brownian motions and independent compound Poisson processes which are supposed to be mutually independent, too. That is,

$$H_t^j = \mu_j t + \beta_j \vartheta_t^j + \sigma_j B_{\vartheta_t^j}^j - Z_t^j, \quad H_0^j = 0, \qquad (1)$$

where $\mu_j, \beta_j \in \mathbb{R}$, $\sigma_j \geq 0$, $(B_t^j)_{t \geq 0}$ are the Wiener processes correlated with coefficients $\rho_{jl}$, $(\vartheta_t^j)_{t \leq T}$ are independent with the Wiener processes subordinators and $Z_t^j = \sum_{l=0}^{N_t^j} \xi_{jl}$, $\xi_{j0} \equiv 0$, where $(N_t^j)_{t \leq T}$, $N_0^j = 0$, are the Poisson processes with intensities $\lambda_j$ and $\xi_{jl} \geq 0$, $l = 1, 2, \ldots$, are independent arbitrary identically distributed for every $j$ random variables, where $j$ is the number of asset. Throughout this paper, the problem of pricing of digital asset-or-nothing call option in foreign currency, namely which has the payoff function

$$\mathbb{DC}_T = S_T^3 S_T^2 I_{\{S_T^1 \geq K\}}, \quad K > 0, \qquad (2)$$

is discussed. The dynamics $S_t^3$ relates here to the exchange rate between the domestic and the foreign currencies. The stock prices $S_t^1$ and $S_t^2$ are measured in the domestic currency. It is supposed that the non-risky assets (bank accounts) in domestic and foreign currencies $R_t^d$ and $R_t^f$, $t \leq T$, have fixed interest rates $r_d, r \geq 0$ and $R_t^d = e^{r_d t}$, $R_t^f = e^{rt}$.

It is easy to observe that the problem of pricing the options with payoffs (2) includes the same problem for digital asset-or-nothing and cash-or-nothing call options with payoffs $S_T^1 I_{\{S_T^1 \geq K\}}$ and $\tilde{K} I_{\{S_T^1 \geq K\}}$, $\tilde{K} > 0$, for the options in foreign currency with payoffs $S_T^3(S_T^1 - K)^+$ and for many other options. Indeed, if we discuss for example the payoffs $S_T^3(S_T^1 - K)^+$, we just suppose in (1)–(2) that $\mu_2 = \beta_2 = \sigma_2 = 0$ and $\xi_{2l} \equiv 0$.

Next, it is suggested in our model that the stock prices satisfy the inequality

$$\mathrm{E}(S_T^3 S_T^2) < \infty. \tag{3}$$

Let

$$X_t^j = \mu_j t + \beta_j \vartheta_t^j + \sigma_j B_{\vartheta_t^j}^j. \tag{4}$$

Then

$$\mathrm{E}\bigl(e^{X_T^2 + X_T^3} | \vartheta_T^2, \vartheta_T^3\bigr) = e^{\sum_{j=2}^3 (\mu_j T + \beta_j \vartheta_T^j) + \frac{\sum_{j=2}^3 \sigma_j^2 \vartheta_T^j + 2\rho_{23}\sigma_2\sigma_3 \sqrt{\vartheta_T^2 \vartheta_T^3}}{2}}$$

and hence (3) is equivalent to

$$\mathrm{E}\left(e^{\sum_{j=2}^3 \beta_j \vartheta_T^j + \frac{\sum_{j=2}^3 \sigma_j^2 \vartheta_T^j + 2\rho_{23}\sigma_2\sigma_3 \sqrt{\vartheta_T^2 \vartheta_T^3}}{2}}\right) < \infty. \tag{5}$$

Since our model is not the classical two-asset financial market model (see for example the book by Shiryaev [41]), we need to consider at first the question of hedging of the option with payoffs (2). There are four hedging instruments in our situation. Namely, the bank account in foreign currency $R_t^f$, the bank account in domestic currency transferred in foreign currency with the dynamics $S_t^3 R_t^d$ and the two stocks in foreign currency $S_t^3 S_t^1$ and $S_t^3 S_t^2$. Leaving aside a well-investigated in literature question of change of measure (see for example Eberlein et al. [12], Kallsen and Shiryaev [22], Madan and Milne [30], Ch. VII.3 of Shiryaev [41] and Ch. 6 of Schoutens [39]), let us assume that the all four assets discounted with respect to the bank account in foreign currency (i.e., the processes $R_t^f/R_t^f \equiv 1$, $S_t^3 R_t^d/R_t^f$, $S_t^3 S_t^1/R_t^f$, $S_t^3 S_t^2/R_t^f$) are martingales with respect to the initial probability measure. Then the price of the option with payoffs (2) is

$$\mathbb{DC} = e^{-rT} \mathrm{E}(\mathbb{DC}_T) = e^{-rT} \mathrm{E}\bigl(S_T^3 S_T^2 I_{\{S_T^1 \geq K\}}\bigr). \tag{6}$$

**Remark 1.** Similarly to (2), the digital asset-or-nothing put option in foreign currency has the payoffs at expiry

$$\mathbb{DP}_T = S_T^3 S_T^2 I_{\{S_T^1 < K\}}, \quad K > 0.$$

Hence its price

$$\mathbb{DP} = e^{-rT} \mathrm{E}(\mathbb{DP}_T) = e^{-rT} \mathrm{E}\bigl(S_T^3 S_T^2 I_{\{S_T^1 < K\}}\bigr)$$
$$= e^{-rT} \mathrm{E}\bigl(S_T^3 S_T^2\bigr) - e^{-rT} \mathrm{E}\bigl(S_T^3 S_T^2 I_{\{S_T^1 \geq K\}}\bigr) = e^{-rT} \mathrm{E}\bigl(S_T^3 S_T^2\bigr) - \mathbb{DC}.$$

For the typical case of put option in foreign currency we have for its price the identity

$$\mathbb{P} = e^{-rT} \mathrm{E}\bigl(S_T^3(K - S_T^1)^+\bigr) = e^{-rT} \mathrm{E}\bigl(S_T^3(K - S_T^1)\bigr) - e^{-rT} \mathrm{E}\bigl(S_T^3(K - S_T^1)^-\bigr)$$

$$= e^{-rT}\mathrm{E}\big(S_T^3(K - S_T^1)\big) + e^{-rT}\mathrm{E}\big(S_T^3(S_T^1 - K)^+\big)$$
$$= e^{-rT}K\mathrm{E}\big(S_T^3\big) - e^{-rT}\mathrm{E}\big(S_T^3 S_T^1\big) + \mathbb{C},$$

where $\mathbb{C}$ is the price of call option in foreign currency. That is, results for the prices of call options in foreign currency can be exploited for the computing of prices of put options as well.

Next, we introduce some necessary notations. We denote as

$$\mathrm{N}(u), \quad u \in \mathbb{R}, \qquad \Gamma(u), \quad u > 0, \qquad \mathrm{B}(u_1, u_2), \quad u_1 > 0, u_2 > 0$$

and

$$\mathrm{M}_{u_1}(u_2), \quad u_1 \in \mathbb{R}, u_2 > 0,$$

the normal distribution function, the gamma function, the beta function and the MacDonald function (the modified Bessel function of the second kind), respectively. The hypergeometric Gauss function is denoted as

$$\mathrm{G}(u_1, u_2, u_3; u_4), \quad u_1, u_2, u_3 \in \mathbb{R}, u_4 < 1.$$

Also, the degenerate Appell function (or the Humbert series) which is the double sum

$$\mathrm{A}(u_1, u_2, u_3; u_4, u_5) = \sum_{m=0}^{\infty}\sum_{n=0}^{\infty} \frac{(u_1)_{m+n}(u_2)_m}{m! n! (u_3)_{m+n}} u_4^m u_5^n$$

with $u_1, u_2, u_3, u_5 \in \mathbb{R}$ and $|u_4| < 1$, where $(u)_l, l \in \mathbb{N} \cup \{0\}$, is the Pochhammer's symbol, is exploited. For more information on the special mathematical functions above, see Bateman and Erdélyi [7], Gradshteyn and Ryzhik [16], Whittaker and Watson [44].

## 3 Gamma time change

The gamma process $\gamma_t = \gamma_t(a, b), a > 0, b > 0$, is a purely discontinuous Lévy process with gamma-distributed increments and $\gamma_0 = 0$. It is the subordinator with the probability density function

$$f(\gamma_t, x) = \frac{b^{at} x^{at-1} e^{-bx}}{\Gamma(at)}, \quad x > 0.$$

The gamma process has mean $at/b$ and variance $at/b^2$. If $u < b$, the moment-generating function of the gamma process is

$$\mathrm{E}e^{u\gamma_t} = \left(\frac{b}{b-u}\right)^{at}. \tag{7}$$

For more properties of this process, see the paper by Yor [45] or the monograph by Applebaum [2].

Throughout this section, we assume that the subordinators in (1) and (4) are the gamma processes with unit mean rate, i.e.

$$\vartheta_t^j = \gamma_t^j(a_j) = \gamma_t^j(a_j, a_j). \tag{8}$$

Then the processes $X_t^j$ in (4) become the variance-gamma processes, see Madan et al. [33] or Seneta [40] for more details.

To model dependencies in the subordinators, let us assume that in (8) the subordinators

$$\gamma_t^j = \kappa_j \gamma_t(a) + \kappa_{j1}\gamma_t^1 + \tilde{\kappa}_j \tilde{\gamma}_t^j(\tilde{a}_j), \quad j = 2, 3, \tag{9}$$

where all the gamma processes with unit mean rate $\gamma_t, \gamma_t^1, \tilde{\gamma}_t^2, \tilde{\gamma}_t^3$ are mutually independent, $\kappa_j, \kappa_{j1}, \tilde{\kappa}_j \geq 0$ and $\kappa_j + \kappa_{j1} + \tilde{\kappa}_j = 1, j = 2, 3$.

Since for a gamma distribution $\gamma$ the identity

$$u\gamma(u_1, u_2) \stackrel{Law}{=} \gamma\left(u_1, \frac{u_2}{u}\right)$$

is satisfied, we have from (9) that $a_j = \frac{a}{\kappa_j}$ if $\kappa_j \neq 0$, $a_j = \frac{a_1}{\kappa_{j1}}$ if $\kappa_{j1} \neq 0$, $a_j = \frac{\tilde{a}_j}{\tilde{\kappa}_j}$ if $\tilde{\kappa}_j \neq 0$ and hence the equality

$$\left(a_j - \frac{a}{\kappa_j}\right)I_{\{\kappa_j > 0\}} = \left(a_j - \frac{a_1}{\kappa_{j1}}\right)I_{\{\kappa_{j1} > 0\}} = \left(a_j - \frac{\tilde{a}_j}{\tilde{\kappa}_j}\right)I_{\{\tilde{\kappa}_j > 0\}} = 0 \tag{10}$$

holds. Next, because the identity

$$\gamma(u_1, u) + \tilde{\gamma}(u_2, u) \stackrel{Law}{=} \gamma(u_1 + u_2, u)$$

holds for arbitrary independent gamma distributions $\gamma$ and $\tilde{\gamma}$, one could observe from (9), as $\gamma_t, \gamma_t^1, \tilde{\gamma}_t^j$ are mutually independent, that

$$a_j = aI_{\{\kappa_j > 0\}} + a_1 I_{\{\kappa_{j1} > 0\}} + \tilde{a}_j I_{\{\tilde{\kappa}_j > 0\}} \tag{11}$$

in our model, $j = 2, 3$. Alternatively, the identities (10) and (11) can be seen from the equality for characteristic functions of (9)

$$\left(\frac{a_j}{a_j - iu}\right)^{a_j t} = \left(\frac{a/\kappa_j}{a/\kappa_j - iu}\right)^{at} \left(\frac{a_1/\kappa_{j1}}{a_1/\kappa_{j1} - iu}\right)^{a_1 t} \left(\frac{\tilde{a}_j/\tilde{\kappa}_j}{\tilde{a}_j/\tilde{\kappa}_j - iu}\right)^{\tilde{a}_j t}$$

if all $\kappa_j > 0$, $\kappa_{j1} > 0$, $\tilde{\kappa}_j > 0$. The theorem below gives us the price (6) in the case of the independent Brownian motions in (1).

**Theorem 1.** *Let the stock log-returns be defined in* (1), *the subordinators $\vartheta_t^j$ be gamma distributed, satisfy* (8)–(9), *and $\rho_{12} = \rho_{13} = \rho_{23} = 0$. Set*

$$b = \sum_{j=2}^{3} \kappa_{j1}\left(\beta_j + \frac{\sigma_j^2}{2}\right).$$

*Then the double inequality for the price* (6)

$$\sum_{n_1=0}^{N_1} \sum_{n_2=0}^{N_2} \sum_{n_3=0}^{N_3} \frac{\lambda_1^{n_1} \lambda_2^{n_2} \lambda_3^{n_3} T^{n_1+n_2+n_3} e^{-(\lambda_1+\lambda_2+\lambda_3)T} \mathrm{DC}(n_1,n_2,n_3)}{n_1! n_2! n_3!}$$
$$\leq \mathbb{DC}$$
$$\leq \sum_{n_1=0}^{N_1} \sum_{n_2=0}^{N_2} \sum_{n_3=0}^{N_3} \frac{\lambda_1^{n_1} \lambda_2^{n_2} \lambda_3^{n_3} T^{n_1+n_2+n_3} e^{-(\lambda_1+\lambda_2+\lambda_3)T} \mathrm{DC}(n_1,n_2,n_3)}{n_1! n_2! n_3!}$$
$$+ \mathrm{DC}(N_1, N_2, N_3) \left(1 - \sum_{n_1=0}^{N_1} \frac{\lambda_1^{n_1} e^{-\lambda_1 T}}{n_1!}\right)$$
$$\times \left(1 - \sum_{n_2=0}^{N_2} \frac{\lambda_2^{n_2} e^{-\lambda_2 T}}{n_2!}\right) \left(1 - \sum_{n_3=0}^{N_3} \frac{\lambda_3^{n_3} e^{-\lambda_3 T}}{n_3!}\right) \qquad (12)$$

*holds for any* $N_1, N_2, N_3$ *with a decreasing function* $\mathrm{DC}(n_1, n_2, n_3)$ *and*

$$\begin{aligned}
&\mathrm{DC}(n_1, n_2, n_3) \\
&= \frac{e^{(\mu_2+\mu_3-r)T} a_1^{a_1 T}}{(a_1-b)^{a_1 T} \Gamma(a_1 T)\sqrt{2\pi}} \mathrm{E}\big(e^{-\sum_{l=0}^{n_2} \xi_{2l}}\big) \mathrm{E}\big(e^{-\sum_{l=0}^{n_3} \xi_{3l}}\big) \left(\frac{\tilde{a}_2}{\tilde{a}_2 - \tilde{\kappa}_2(\beta_2 + \frac{\sigma_2^2}{2})}\right)^{\tilde{a}_2 T} \\
&\quad \times \left(\frac{\tilde{a}_3}{\tilde{a}_3 - \tilde{\kappa}_3(\beta_3 + \frac{\sigma_3^2}{2})}\right)^{\tilde{a}_3 T} \left(\frac{a}{a - \sum_{j=2}^{3} \kappa_j(\beta_j + \frac{\sigma_j^2}{2})}\right)^{aT} \\
&\quad \times \left(\Lambda \mathrm{P}\left(\sum_{l=0}^{n_1} \xi_{1l} = \mu_1 T - K\right) + \mathrm{E}\left(\Xi\left(\sum_{l=0}^{n_1} \xi_{1l}\right) I_{\{\sum_{l=0}^{n_1} \xi_{1l} \neq \mu_1 T - K\}}\right)\right),
\end{aligned}$$

*where*

$$\begin{aligned}
\Lambda = \Gamma\left(a_1 T + \frac{1}{2}\right) &\left(\frac{\mathrm{B}(\frac{1}{2}, a_1 T)}{\sqrt{2}}\right. \\
&\left. + \frac{\beta_1}{\sigma_1 \sqrt{a_1 - b}} \mathrm{G}\left(a_1 T + \frac{1}{2}, \frac{1}{2}, \frac{3}{2}; -\frac{\beta_1^2}{2(a_1-b)\sigma_1^2}\right)\right)
\end{aligned} \qquad (13)$$

*and*

$$\Xi(x) = |s|^{a_1 T - \frac{1}{2}} e^s (1+q)^{a_1 T} \big(\mathrm{B}(a_1 T, 1)\big(|s| \mathrm{M}_{a_1 T + \frac{1}{2}}(|s|) \\
+ s \mathrm{M}_{a_1 T - \frac{1}{2}}(|s|)\big) \mathrm{A}_0 - (1+q) s \mathrm{B}(a_1 T + 1, 1) \mathrm{M}_{a_1 T - \frac{1}{2}}(|s|) \mathrm{A}_1 \big) \qquad (14)$$

*with*

$$q = \frac{\beta_1}{\sqrt{\beta_1^2 + 2(a_1-b)\sigma_1^2}}, \qquad s = s(x) = \frac{(\mu_1 T - K - x)\sqrt{\beta_1^2 + 2(a_1-b)\sigma_1^2}}{\sigma_1}$$

*and* $\mathrm{A}_j = \mathrm{A}(a_1 T + j, 1 - a_1 T, a_1 T + 1 + j; \frac{1+q}{2}, -s(1+q))$.

The following example illustrates how Theorem 1 works when $Z_t^j$ are standard Poisson processes.

**Example 1.** Let $\xi_{jl} \equiv \varpi_j$, $j = 1, 2, 3$, $l = 1, 2, \ldots$, where $\varpi_j \geq 0$ are constants. Then $Z_t^j \equiv \varpi_j N_t^j$ (Poisson processes) and the result of Theorem 1 holds with

$$\mathrm{DC}(n_1, n_2, n_3) = \frac{e^{(\mu_2+\mu_3-r)T - \varpi_2 n_2 - \varpi_3 n_3} a_1^{a_1 T}}{(a_1 - b)^{a_1 T} \Gamma(a_1 T)\sqrt{2\pi}} \left(\frac{\tilde{a}_2}{\tilde{a}_2 - \tilde{\kappa}_2(\beta_2 + \frac{\sigma_2^2}{2})}\right)^{\tilde{a}_2 T}$$

$$\times \left(\frac{\tilde{a}_3}{\tilde{a}_3 - \tilde{\kappa}_3(\beta_3 + \frac{\sigma_3^2}{2})}\right)^{\tilde{a}_3 T} \left(\frac{a}{a - \sum_{j=2}^{3} \kappa_j(\beta_j + \frac{\sigma_j^2}{2})}\right)^{aT}$$

$$\times \left(\Lambda I_{\{\varpi_1 n_1 = \mu_1 T - K\}} + \Xi(\varpi_1 n_1) I_{\{\varpi_1 n_1 \neq \mu_1 T - K\}}\right).$$

Theorem 2 computes us the price (6) in the case when the exchange rate $S_t^3$ and the underlying asset $S_t^2$ are strongly dependent but the indicator stock $S_t^1$ is weakly dependent on them.

**Theorem 2.** *Assume that in* (1) $\rho_{12} = \rho_{13} = 0$, *the subordinators are gamma distributed, satisfy* (8)–(9), *and* $\gamma_t^3 = \gamma_t^2 = \kappa_2 \gamma_t + \kappa_{21} \gamma_t^1$. *Let*

$$b = \kappa_{21} \left[\sum_{j=2}^{3} \left(\beta_j + \frac{\sigma_j^2}{2}\right) + \rho_{23} \sigma_2 \sigma_3\right].$$

*Then* (12) *is satisfied with*

$$\mathrm{DC}(n_1, n_2, n_3) = \frac{e^{(\mu_2+\mu_3-r)T} a_1^{a_1 T}}{(a_1 - b)^{a_1 T} \Gamma(a_1 T)\sqrt{2\pi}} \mathrm{E}\left(e^{-\sum_{l=0}^{n_2} \xi_{2l}}\right) \mathrm{E}\left(e^{-\sum_{l=0}^{n_3} \xi_{3l}}\right)$$

$$\times \left(\frac{a}{a - \kappa_2\left[\sum_{j=2}^{3}(\beta_j + \frac{\sigma_j^2}{2}) + \rho_{23}\sigma_2\sigma_3\right]}\right)^{aT}$$

$$\times \left(\Lambda \mathrm{P}\left(\sum_{l=0}^{n_1} \xi_{1l} = \mu_1 T - K\right)\right.$$

$$\left.+ \mathrm{E}\left(\Xi\left(\sum_{l=0}^{n_1} \xi_{1l}\right) I_{\{\sum_{l=0}^{n_1} \xi_{1l} \neq \mu_1 T - K\}}\right)\right),$$

*where $\Lambda$ and $\Xi(x)$ are defined in* (13) *and* (14), *respectively.*

The next theorem considers the case when all risky assets are strongly dependent.

**Theorem 3.** *Let the subordinators in* (1) *be gamma distributed, satisfy* (8)–(9), *and* $\gamma_t^3 = \gamma_t^2 = \gamma_t^1$. *Set*

$$b = \sum_{j=2}^{3} \left(\beta_j + \frac{\sigma_j^2}{2}\right) + \rho_{23} \sigma_2 \sigma_3.$$

*Then* (12) *holds with*

$$\mathrm{DC}(n_1, n_2, n_3) = \frac{e^{(\mu_2+\mu_3-r)T} a_1^{a_1 T}}{(a_1-b)^{a_1 T}\Gamma(a_1 T)\sqrt{2\pi}} \mathrm{E}\bigl(e^{-\sum_{l=0}^{n_2}\xi_{2l}}\bigr)\mathrm{E}\bigl(e^{-\sum_{l=0}^{n_3}\xi_{3l}}\bigr)$$

$$\times \Biggl(\Lambda \mathrm{P}\Biggl(\sum_{l=0}^{n_1}\xi_{1l} = \mu_1 T - K\Biggr)$$

$$+ \mathrm{E}\Biggl(\Xi\Biggl(\sum_{l=0}^{n_1}\xi_{1l}\Biggr) I_{\{\sum_{l=0}^{n_1}\xi_{1l}\neq \mu_1 T-K\}}\Biggr)\Biggr),$$

*where*

$$\Lambda = \Gamma\Bigl(a_1 T + \frac{1}{2}\Bigr)\Biggl(\frac{\mathrm{B}(\frac{1}{2}, a_1 T)}{\sqrt{2}}$$
$$+ \frac{\beta_1 + \sum_{j=2}^{3}\rho_{1j}\sigma_1\sigma_j}{\sigma_1\sqrt{a_1-b}}\mathrm{G}\Bigl(a_1 T + \frac{1}{2}, \frac{1}{2}, \frac{3}{2}; -\frac{(\beta_1 + \sum_{j=2}^{3}\rho_{1j}\sigma_1\sigma_j)^2}{2(a_1-b)\sigma_1^2}\Bigr)\Biggr) \quad (15)$$

*and* $\Xi(x)$ *is defined by* (14) *with*

$$q = \frac{\beta_1 + \sum_{j=2}^{3}\rho_{1j}\sigma_1\sigma_j}{\sqrt{(\beta_1 + \sum_{j=2}^{3}\rho_{1j}\sigma_1\sigma_j)^2 + 2(a_1-b)\sigma_1^2}}$$

*and*

$$s = s(x) = \frac{(\mu_1 T - K - x)\sqrt{(\beta_1 + \sum_{j=2}^{3}\rho_{1j}\sigma_1\sigma_j)^2 + 2(a_1-b)\sigma_1^2}}{\sigma_1}.$$

Example 2 shows how Theorem 3 can be applied to the problem of pricing of the standard European call option in foreign currency which has the payoffs at expiry $S_T^2(S_T^1 - K)^+$.

**Example 2.** Assume that $S_t^3 \equiv S_t^1$ and $\xi_{jl} \equiv 0$, $j = 1, 2, 3$, $l = 1, 2, \ldots$ under the conditions of Theorem 3. Then

$$\mathbb{DC} = \mathrm{DC}(0,0,0) = \frac{e^{(\mu_2+\mu_1-r)T}a_1^{a_1 T}}{(a_1-b)^{a_1 T}\Gamma(a_1 T)\sqrt{2\pi}}(\Lambda I_{\{\mu_1 T=K\}} + \Xi I_{\{\mu_1 T\neq K\}}), \quad (16)$$

*where*

$$\Lambda = \Gamma\Bigl(a_1 T + \frac{1}{2}\Bigr)\Biggl(\frac{\mathrm{B}(\frac{1}{2}, a_1 T)}{\sqrt{2}}$$
$$+ \frac{\beta_1 + \sigma_1^2 + \rho_{12}\sigma_1\sigma_2}{\sigma_1\sqrt{a_1-b}}\mathrm{G}\Bigl(a_1 T + \frac{1}{2}, \frac{1}{2}, \frac{3}{2}; -\frac{(\beta_1 + \sigma_1^2 + \rho_{12}\sigma_1\sigma_2)^2}{2(a_1-b)\sigma_1^2}\Bigr)\Biggr)$$

with

$$b = \sum_{j=1}^{2}\left(\beta_j + \frac{\sigma_j^2}{2}\right) + \rho_{12}\sigma_1\sigma_2$$

and $\Xi$ is set by (14) with

$$q = \frac{\beta_1 + \sigma_1^2 + \rho_{12}\sigma_1\sigma_2}{\sqrt{(\beta_1 + \sigma_1^2 + \rho_{12}\sigma_1\sigma_2)^2 + 2(a_1 - b)\sigma_1^2}}$$

and

$$s = \frac{(\mu_1 T - K)\sqrt{(\beta_1 + \sigma_1^2 + \rho_{12}\sigma_1\sigma_2)^2 + 2(a_1 - b)\sigma_1^2}}{\sigma_1}.$$

Next, let $S_t^3 \equiv 1$, $\xi_{jl} \equiv 0$, $j = 1, 2$, $l = 1, 2, \ldots$ and the conditions of Theorem 3 hold. Then

$$\mathbb{DC} = \mathrm{DC}(0,0) = \frac{e^{(\mu_2-r)T}a_1^{a_1 T}}{(a_1-b)^{a_1 T}\Gamma(a_1 T)\sqrt{2\pi}}(\Lambda I_{\{\mu_1 T = K\}} + \Xi I_{\{\mu_1 T \neq K\}}), \quad (17)$$

where

$$\Lambda = \Gamma\left(a_1 T + \frac{1}{2}\right)\left(\frac{\mathrm{B}(\frac{1}{2}, a_1 T)}{\sqrt{2}}\right.$$
$$\left. + \frac{\beta_1 + \rho_{12}\sigma_1\sigma_2}{\sigma_1\sqrt{a_1 - b}}\mathrm{G}\left(a_1 T + \frac{1}{2}, \frac{1}{2}, \frac{3}{2}; -\frac{(\beta_1 + \rho_{12}\sigma_1\sigma_2)^2}{2(a_1 - b)\sigma_1^2}\right)\right)$$

with $b = \beta_2 + \frac{\sigma_2^2}{2}$ and $\Xi$ is defined in (14) with

$$q = \frac{\beta_1 + \rho_{12}\sigma_1\sigma_2}{\sqrt{(\beta_1 + \rho_{12}\sigma_1\sigma_2)^2 + 2(a_1 - b)\sigma_1^2}}$$

and

$$s = \frac{(\mu_1 T - K)\sqrt{(\beta_1 + \rho_{12}\sigma_1\sigma_2)^2 + 2(a_1 - b)\sigma_1^2}}{\sigma_1}.$$

Combining together (16) and (17), one can obtain the result of Theorem 1 from Ivanov and Ano [20].

Now we will consider the case when the indicator stock $S_t^1$ and the exchange rate $S_t^3$ are strongly dependent but the underlying asset $S_t^2$ is weakly dependent on them.

**Theorem 4.** *Assume that in* (1) $\rho_{23} = \rho_{12} = 0$, *the subordinators are gamma distributed, satisfy* (8)–(9), *and* $\gamma_t^3 = \gamma_t^1$, $\gamma_t^2 = \kappa_{21}\gamma_t^1 + \tilde{\kappa}_2\tilde{\gamma}_t^2$. *Let*

$$b = \beta_3 + \frac{\sigma_3^2}{2} + \kappa_{21}\left(\beta_2 + \frac{\sigma_2^2}{2}\right).$$

*Then* (12) *is satisfied with*

$$\mathrm{DC}(n_1, n_2, n_3) = \frac{e^{(\mu_2+\mu_3-r)T} a_1^{a_1 T}}{(a_1-b)^{a_1 T} \Gamma(a_1 T)\sqrt{2\pi}} \mathrm{E}\big(e^{-\sum_{l=0}^{n_2} \xi_{2l}}\big)\mathrm{E}\big(e^{-\sum_{l=0}^{n_3} \xi_{3l}}\big)$$

$$\times \left(\frac{\tilde{a}_2}{\tilde{a}_2 - \tilde{\kappa}_2(\beta_2 + \frac{\sigma_2^2}{2})}\right)^{\tilde{a}_2 T} \left(\Lambda \mathrm{P}\left(\sum_{l=0}^{n_1} \xi_{1l} = \mu_1 T - K\right)\right.$$

$$\left. + \mathrm{E}\left(\Xi\left(\sum_{l=0}^{n_1} \xi_{1l}\right) I_{\{\sum_{l=0}^{n_1} \xi_{1l} \neq \mu_1 T - K\}}\right)\right),$$

*where*

$$\Lambda = \Gamma\left(a_1 T + \frac{1}{2}\right)\left(\frac{\mathrm{B}(\frac{1}{2}, a_1 T)}{\sqrt{2}}\right.$$
$$\left. + \frac{\beta_1 + \sigma_1\sigma_3}{\sigma_1\sqrt{a_1 - b}} \mathrm{G}\left(a_1 T + \frac{1}{2}, \frac{1}{2}, \frac{3}{2}; -\frac{(\beta_1 + \sigma_1\sigma_3)^2}{2(a_1 - b)\sigma_1^2}\right)\right)$$

*and $\Xi(x)$ is defined by* (14) *with*

$$q = \frac{\beta_1 + \sigma_1\sigma_3}{\sqrt{(\beta_1 + \sigma_1\sigma_3)^2 + 2(a_1 - b)\sigma_1^2}}$$

*and*

$$s = s(x) = \frac{(\mu_1 T - K - x)\sqrt{(\beta_1 + \sigma_1\sigma_3)^2 + 2(a_1 - b)\sigma_1^2}}{\sigma_1}.$$

**Remark 2.** One could notice that the result symmetric to Theorem 4 can be established. It should be assumed then that the indicator stock $S_t^1$ and the underlying asset $S_t^2$ are strongly dependent but the exchange rate $S_t^3$ is weakly dependent on them. That is, the conditions $\rho_{23} = \rho_{13} = 0$ and $\gamma_t^2 = \gamma_t^1$, $\gamma_t^3 = \kappa_{31}\gamma_t^1 + \tilde{\kappa}_3\tilde{\gamma}_t^3$ have to be proposed.

## 4 Inverse-Gaussian time change

Let $(\tilde{B}_s)_{s \geq 0}$ be a Brownian motion, $\phi > 0$ and $a \geq 0$. Set for $t \geq 0$

$$\varkappa_t = \varkappa_t(\phi, a) = \inf\{s \geq 0 : \tilde{B}_s + as \geq \phi t\}. \tag{18}$$

The subordinator $(\varkappa_t)_{t \geq 0}$ is called the inverse-Gaussian process and has the probability density function

$$f(\varkappa_t, x) = \frac{\phi t}{\sqrt{2\pi}} x^{-\frac{3}{2}} e^{a\phi t - \frac{1}{2}(a^2 x + \frac{(\phi t)^2}{x})}, \tag{19}$$

see, for example, (1.26) in Applebaum [2]. The mean of $\varkappa_t$ is

$$\mathrm{E}(\varkappa_t) = \frac{\phi t}{\sqrt{2\pi}} e^{a\phi t} \int_0^\infty \sqrt{x} e^{-\frac{1}{2}(a^2 x + \frac{(\phi t)^2}{x})} dx$$

$$= e^{a\phi t} \frac{(\phi t)^{\frac{3}{2}}}{\sqrt{a}} \sqrt{\frac{2}{\pi}} \mathrm{M}_{\frac{1}{2}}(a\phi t) = \frac{\phi t}{a}$$

with respect to 3.471.9 and 8.469.3 from Gradshteyn and Ryzhik [16]. In this section we assume that the subordinator in (1) and (4) is the inverse-Gaussian process with unit mean rate, that is, we set

$$\vartheta_t^j = \varkappa_t^j(\phi_j) = \varkappa_t^j(\phi_j, \phi_j). \tag{20}$$

Then the processes $X_t^j$ in (4) become the normal-inverse Gaussian processes, see, for example, Ivanov and Temnov [21] and references therein or Applebaum [2].

Similarly to (9), we assume that

$$\varkappa_t^j = \kappa_j \varkappa_t(\phi) + \kappa_{j1} \varkappa_t^1 + \tilde{\kappa}_j \tilde{\varkappa}_t^j(\tilde{\phi}_j), \quad j = 2, 3, \tag{21}$$

where all the inverse-Gaussian processes with unit mean rate $\varkappa_t, \varkappa_t^1, \tilde{\varkappa}_t^2, \tilde{\varkappa}_t^3$ are mutually independent, $\kappa_j, \kappa_{j1}, \tilde{\kappa}_j \geq 0$ and $\kappa_j + \kappa_{j1} + \tilde{\kappa}_j = 1$, $j = 2, 3$. Because for arbitrary independent inverse-Gaussian distributions $\varkappa$ and $\tilde{\varkappa}$ the identities

$$u\varkappa(u_1, u_2) \stackrel{Law}{=} \varkappa\left(u_1 \sqrt{u}, \frac{u_2}{\sqrt{u}}\right)$$

and

$$\varkappa(u_1, u) + \tilde{\varkappa}(u_2, u) \stackrel{Law}{=} \varkappa(u_1 + u_2, u)$$

are satisfied, one could observe that in the model (21)

$$\left(\phi_j - \frac{\phi}{\sqrt{\kappa_j}}\right) I_{\{\kappa_j > 0\}} = \left(\phi_j - \frac{\phi_1}{\sqrt{\kappa_{j1}}}\right) I_{\{\kappa_{j1} > 0\}}$$

$$= \left(\phi_j - \frac{\tilde{\phi}_j}{\sqrt{\tilde{\kappa}_j}}\right) I_{\{\tilde{\kappa}_j > 0\}} = 0$$

and

$$\phi_j = \phi\sqrt{\kappa_j} I_{\{\kappa_j > 0\}} + \phi_1 \sqrt{\kappa_{j1}} I_{\{\kappa_{j1} > 0\}} + \tilde{\phi}_j \sqrt{\tilde{\kappa}_j} I_{\{\tilde{\kappa}_j > 0\}}.$$

The next theorem suggests the conditions of dependence in (1) which are similar to those of Theorem 1.

**Theorem 5.** *Let in* (1) $\rho_{12} = \rho_{13} = \rho_{23} = 0$, *the subordinators satisfy* (20)–(21), *and*

$$\phi_1^2 = 2 \sum_{j=2}^3 \kappa_{j1} \left(\beta_j + \frac{\sigma_j^2}{2}\right).$$

*Then* (12) *holds with*

$$\mathrm{DC}(n_1, n_2, n_3) = \frac{e^{(\mu_2+\mu_3-r+\phi_1^2)T}}{2\sqrt{\pi}} \mathrm{E}\bigl(e^{-\sum_{l=0}^{n_2}\xi_{2l}}\bigr)\mathrm{E}\bigl(e^{-\sum_{l=0}^{n_3}\xi_{3l}}\bigr)$$

$$\times e^{T\bigl(\tilde{\phi}_2(\tilde{\phi}_2-\sqrt{\tilde{\phi}_2^2-2\tilde{\kappa}_2(\beta_2+\frac{\sigma_2^2}{2})})+\tilde{\phi}_3(\tilde{\phi}_3-\sqrt{\tilde{\phi}_3^2-2\tilde{\kappa}_3(\beta_3+\frac{\sigma_3^2}{2})})\bigr)}$$

$$\times e^{\phi T\bigl(\phi-\sqrt{\phi^2-2\sum_{j=2}^{3}\kappa_j(\beta_j+\frac{\sigma_j^2}{2})}\bigr)}$$

$$\times \left(\mathrm{E}\Lambda\left(\sum_{l=0}^{n_1}\xi_{1l}\right)I_{\{\beta_1=0\}} + \mathrm{E}\Xi\left(\sum_{l=0}^{n_1}\xi_{1l}\right)I_{\{\beta_1\neq 0\}}\right),$$

*where*

$$\Lambda(x) = \sqrt{\pi} + \frac{2}{\sqrt{\pi}}\mathrm{sign}(\mu_1 T - K - x)\arctan\left(\frac{|\mu_1 T - K - x|}{\sigma_1\phi_1 T}\right) \qquad (22)$$

*and*

$$\Xi(x) = \frac{|\varsigma|e^{|\varsigma|}}{\sqrt{q+1}}\bigl(\mathrm{M}_1(|\varsigma|)\Upsilon_0 + \mathrm{M}_0(|\varsigma|)(\Upsilon_0 - (q+1)\Upsilon_1)\bigr) \qquad (23)$$

*with*

$$\varsigma = \varsigma(x) = \frac{\beta_1}{\sigma_1^2}\sqrt{(\mu_1 T - K - x)^2 + (\sigma_1\phi_1 T)^2},$$

$$q = q(x) = \frac{\mu_1 T - K - x}{\sqrt{(\mu_1 T - K - x)^2 + (\sigma_1\phi_1 T)^2}} \quad \text{and}$$

$$\Upsilon_j = \Upsilon_j(x) = \mathrm{B}\left(\frac{1}{2}+j, 1\right)\mathrm{A}\left(\frac{1}{2}+j, \frac{1}{2}, \frac{3}{2}+j; \frac{q+1}{2}, -|\varsigma|(q+1)\right).$$

The following example applies the result of Theorem 5 to the case of standard Poisson processes.

**Example 3.** Assume that $\xi_{jl} \equiv \varpi_j$, $j = 1, 2, 3$, $l = 1, 2, \ldots$, where $\varpi_j \geq 0$ are constants. Then $Z_t^j \equiv \varpi_j N_t^j$ (Poisson processes) and the result of Theorem 1 holds with

$$\mathrm{DC}(n_1, n_2, n_3) = \frac{e^{(\mu_2+\mu_3-r+\phi_1^2)T-\varpi_2 n_2-\varpi_3 n_3}}{2\sqrt{\pi}}$$

$$\times e^{T\bigl(\tilde{\phi}_2(\tilde{\phi}_2-\sqrt{\tilde{\phi}_2^2-2\tilde{\kappa}_2(\beta_2+\frac{\sigma_2^2}{2})})+\tilde{\phi}_3(\tilde{\phi}_3-\sqrt{\tilde{\phi}_3^2-2\tilde{\kappa}_3(\beta_3+\frac{\sigma_3^2}{2})})\bigr)}$$

$$\times e^{\phi T\bigl(\phi-\sqrt{\phi^2-2\sum_{j=2}^{3}\kappa_j(\beta_j+\frac{\sigma_j^2}{2})}\bigr)}$$

$$\times \bigl(\Lambda(\varpi_1 n_1)I_{\{\beta_1=0\}} + \Xi(\varpi_1 n_1)I_{\{\beta_1\neq 0\}}\bigr).$$

The next two theorems are analogues of Theorem 2 and Theorem 3, respectively.

**Theorem 6.** *Assume that in* (1) $\rho_{12} = \rho_{13} = 0$, *the subordinators* $\varkappa_t^3 = \varkappa_t^2 = \kappa_2 \varkappa_t + \kappa_{21} \varkappa_t^1$, *and the identity*

$$\phi_1^2 = 2\kappa_{21} \left( \sum_{j=2}^{3} \left( \beta_j + \frac{\sigma_j^2}{2} \right) + \rho_{23} \sigma_2 \sigma_3 \right)$$

*holds for their parameters. Then* (12) *is satisfied with*

$$\mathrm{DC}(n_1, n_2, n_3) = \frac{e^{(\mu_2 + \mu_3 - r + \phi_1^2)T}}{2\sqrt{\pi}} \mathrm{E}\big(e^{-\sum_{l=0}^{n_2} \xi_{2l}}\big) \mathrm{E}\big(e^{-\sum_{l=0}^{n_3} \xi_{3l}}\big)$$

$$\times e^{\phi T \left( \phi - \sqrt{\phi^2 - 2\kappa_2 (\sum_{j=2}^{3}(\beta_j + \frac{\sigma_j^2}{2}) + \rho_{23}\sigma_2\sigma_3)} \right)}$$

$$\times \left( \mathrm{E}\Lambda\left(\sum_{l=0}^{n_1}\xi_{1l}\right) I_{\{\beta_1=0\}} + \mathrm{E}\Xi\left(\sum_{l=0}^{n_1}\xi_{1l}\right) I_{\{\beta_1 \neq 0\}} \right),$$

*where* $\Lambda(x)$ *and* $\Xi(x)$ *are defined in* (22) *and* (23), *respectively*.

**Theorem 7.** *Let the subordinators in* (1) *satisfy*

$$\varkappa_t^3 = \varkappa_t^2 = \varkappa_t^1, \qquad \phi_1^2 = 2\sum_{j=2}^{3}\left(\beta_j + \frac{\sigma_j^2}{2}\right) + \rho_{23}\sigma_2\sigma_3.$$

*Then* (12) *holds with*

$$\mathrm{DC}(n_1, n_2, n_3) = \frac{e^{(\mu_2 + \mu_3 - r + \phi_1^2)T}}{2\sqrt{\pi}} \mathrm{E}\big(e^{-\sum_{l=0}^{n_2}\xi_{2l}}\big) \mathrm{E}\big(e^{-\sum_{l=0}^{n_3}\xi_{3l}}\big)$$

$$\times \left( \mathrm{E}\Lambda\left(\sum_{l=0}^{n_1}\xi_{1l}\right) I_{\{\beta_1 + \rho_{12}\sigma_1\sigma_2 + \rho_{13}\sigma_1\sigma_3 = 0\}} \right.$$

$$\left. + \mathrm{E}\Xi\left(\sum_{l=0}^{n_1}\xi_{1l}\right) I_{\{\beta_1 + \rho_{12}\sigma_1\sigma_2 + \rho_{13}\sigma_1\sigma_3 \neq 0\}} \right),$$

*where* $\Lambda(x)$ *and* $\Xi(x)$ *are defined in* (22) *and* (23) *with*

$$\varsigma = \varsigma(x) = \frac{\beta_1 + \rho_{12}\sigma_1\sigma_2 + \rho_{13}\sigma_1\sigma_3}{\sigma_1^2}\sqrt{(\mu_1 T - K - x)^2 + (\sigma_1\phi_1 T)^2}$$

*and*

$$q = q(x) = \frac{\mu_1 T - K - x}{\sqrt{(\mu_1 T - K - x)^2 + (\sigma_1 \phi_1 T)^2}}.$$

The example below gives us the price of the standard asset-or-nothing digital option computed in Ivanov and Temnov [21] as a corollary of Theorem 7.

**Example 4.** Let $S_t^3 \equiv 1$, $S_t^2 \equiv S_t^1$ and $\xi_{jl} \equiv 0$, $j = 1,2,3$, $l = 1,2,\ldots$. Then $\rho_{12} = 1$, the conditions of Theorem 7 has the form $\phi_1^2 = 2\beta_1 + \sigma_1^2$, $\beta_1 \neq \sigma_1^2$ as in Corollary 3.1 of Ivanov and Temnov [21] and

$$\mathbb{DC} = \mathrm{DC}(0,0) = \frac{e^{(\mu_2+\mu_3-r+\phi_1^2)T}\Xi}{2\sqrt{\pi}},$$

where $\Xi$ is defined in (23) with

$$\varsigma = \frac{\beta_1+\sigma_1^2}{\sigma_1^2}\sqrt{(\mu_1 T - K)^2 + (\sigma_1\phi_1 T)^2} \quad \text{and} \quad q = \frac{\mu_1 T - K}{\sqrt{(\mu_1 T - K)^2 + (\sigma_1\phi_1 T)^2}}.$$

Theorem 8 implies the similar conditions on the dependence between risky assets as Theorem 4 does.

**Theorem 8.** *Assume that $\rho_{23} = \rho_{12} = 0$, $\varkappa_t^3 = \varkappa_t^1$, $\varkappa_t^2 = \kappa_{21}\varkappa_t^1 + \tilde{\kappa}_2\tilde{\varkappa}_t^2$, $\phi_1^2 = 2\beta_3 + \sigma_3^2 + \kappa_{21}(2\beta_2 + \sigma_2^2)$ in (1) under (20)–(21). Then (12) is satisfied with*

$$\mathrm{DC}(n_1, n_2, n_3) = \frac{e^{(\mu_2+\mu_3-r+\phi_1^2)T}}{2\sqrt{\pi}} \mathrm{E}\left(e^{-\sum_{l=0}^{n_2}\xi_{2l}}\right)\mathrm{E}\left(e^{-\sum_{l=0}^{n_3}\xi_{3l}}\right)$$
$$\times e^{\tilde{\phi}_2 T\left(\tilde{\phi}_2 - \sqrt{\tilde{\phi}_2^2 - \tilde{\kappa}_2(2\beta_2+\sigma_2^2)}\right)}$$
$$\times \left(\mathrm{E}\Lambda\left(\sum_{l=0}^{n_1}\xi_{1l}\right)I_{\{\beta_1+\sigma_1\sigma_3=0\}} + \mathrm{E}\Xi\left(\sum_{l=0}^{n_1}\xi_{1l}\right)I_{\{\beta_1+\sigma_1\sigma_3\neq 0\}}\right),$$

*where $\Lambda(x)$ and $\Xi(x)$ are defined in (22) and (23) with*

$$\varsigma = \varsigma(x) = \frac{\beta_1+\sigma_1\sigma_3}{\sigma_1^2}\sqrt{(\mu_1 T - K - x)^2 + (\sigma_1\phi_1 T)^2}$$

*and*

$$q = q(x) = \frac{\mu_1 T - K - x}{\sqrt{(\mu_1 T - K - x)^2 + (\sigma_1\phi_1 T)^2}}.$$

## 5 Conclusion

The paper suggests a foundation for computing of European-style options in the variance-gamma and normal inverse-Gaussian models with extra compound Poisson negative jumps. It is intended to calculate the option prices basing on the knowledge of the price of the digital asset-or-nothing call option in foreign currency. The payoffs of the discussed option build on the values of three risky assets which are assumed to be dependent on each other. Various types of the dependencies between the risky asset prices are considered. The price of the option exploits the values of some special mathematical functions including the hypergeometric ones. A future investigation can relate to discussion of specific types of the compound Poisson process or possibility of the jump in the linear drift, see Ivanov [19].

## 6 Proofs

**Proof of Theorem 1.** We have that the conditional expectation

$$\mathrm{E}\bigl(e^{-rT}S_T^3 S_T^2 I_{\{S_T^1 \geq K\}} | \gamma_T^1, Z_T^1, \gamma_T^2, Z_T^2, \gamma_T^3, Z_T^3\bigr)$$

$$= e^{(\mu_2+\mu_3-r)T - Z_T^2 - Z_T^3 + \sum_{j=1}^{2}\beta_i \gamma_T^j}$$

$$\times \mathrm{E}\left(e^{\sum_{j=2}^{3}\sigma_i\sqrt{\frac{\gamma_T^j}{\gamma_T^1}}B_{\gamma_T^1}^j} I_{\{\mu_1 T + \beta_1 \gamma_T^1 + \sigma_1 B_{\gamma_T^1}^1 - Z_T^1 \geq K\}} | \gamma_T^1, Z_T^1, \gamma_T^2, Z_T^2, \gamma_T^3, Z_T^3\right)$$

$$= e^{(\mu_2+\mu_3-r)T - Z_T^2 - Z_T^3 + \sum_{j=2}^{3}(\beta_j + \frac{\sigma_j^2}{2})\gamma_T^j + \rho_{23}\sigma_2\sigma_3\sqrt{\gamma_T^2 \gamma_T^3}}$$

$$\times \mathrm{E}\left(e^{\sum_{j=2}^{3}\left(\sigma_j\sqrt{\frac{\gamma_T^j}{\gamma_T^1}}B_{\gamma_T^1}^j - \frac{\sigma_j^2 \gamma_T^j}{2}\right) - \rho_{23}\sigma_2\sigma_3\sqrt{\gamma_T^2 \gamma_T^3}}\right.$$

$$\left. \times I_{\{\mu_1 T + \beta_1 \gamma_T^1 + \sigma_1 B_{\gamma_T^1}^1 - Z_T^1 \geq K\}} | \gamma_T^1, Z_T^1, \gamma_T^2, Z_T^2, \gamma_T^3, Z_T^3\right). \tag{24}$$

Let $Q$ be the historical probability measure on the probability space which is generated by the Brownian motions $B_t^j$, $j=1,2,3$, $t \geq 0$. We define a new probability measure $\tilde{Q}$ for fixed trajectories $\gamma_t^j$, $t \leq T$, by the density

$$\frac{d\tilde{Q}_{\gamma_T^1}}{dQ_{\gamma_T^1}} = e^{\sum_{j=2}^{3}\left(\sigma_i\sqrt{\frac{\gamma_T^j}{\gamma_T^1}}B_{\gamma_T^1}^j - \frac{\sigma_j^2 \gamma_T^j}{2}\right) - \rho_{23}\sigma_2\sigma_3\sqrt{\gamma_T^2 \gamma_T^3}}. \tag{25}$$

Then using Corollary 4.5 of [12] one can get that for any $u \in \mathbb{R}$

$$\tilde{Q}\bigl(\log S_T^1 \leq u | \gamma_T^1, Z_T^1, \gamma_T^2, Z_T^2, \gamma_T^3, Z_T^3\bigr)$$

$$= \tilde{Q}\left(\mu_1 T + \left(\beta_1 + \rho_{12}\sigma_1\sigma_2\sqrt{\frac{\gamma_T^2}{\gamma_T^1}} + \rho_{13}\sigma_1\sigma_3\sqrt{\frac{\gamma_T^3}{\gamma_T^1}}\right)\gamma_T^1 + \sigma_1 B_{\gamma_T^1}^{\tilde{Q}} - Z_T^1 \right.$$

$$\left. \leq u | \gamma_T^1, Z_T^1, \gamma_T^2, Z_T^2, \gamma_T^3, Z_T^3\right), \tag{26}$$

where $B_t^{\tilde{Q}}$, $t \leq \gamma_T^1$, is the standard Brownian motion with respect to measure $\tilde{Q}$.

Set

$$\beta_I = \beta_1 + \rho_{12}\sigma_1\sigma_2\sqrt{\frac{\gamma_T^2}{\gamma_T^1}} + \rho_{13}\sigma_1\sigma_3\sqrt{\frac{\gamma_T^3}{\gamma_T^1}}.$$

Then we have from (24) and (26) that

$$\mathrm{E}\bigl(e^{-rT}S_T^3 S_T^2 I_{\{S_T^1 \geq K\}} | \gamma_T^1, Z_T^1, \gamma_T^2, Z_T^2, \gamma_T^3, Z_T^3\bigr)$$

$$= e^{(\mu_2+\mu_3-r)T - Z_T^2 - Z_T^3 + \sum_{j=2}^{3}(\beta_j + \frac{\sigma_j^2}{2})\gamma_T^j + \rho_{23}\sigma_2\sigma_3\sqrt{\gamma_T^2 \gamma_T^3}}$$

$$\times \tilde{Q}\big(\beta_I \gamma_T^1 + \sigma_1 B_{\gamma_T^1}^{\tilde{Q}} \geq K - \mu_1 T + Z_T^1 | \gamma_T^1, Z_T^1, \gamma_T^2, Z_T^2, \gamma_T^3, Z_T^3\big)$$

$$= e^{(\mu_2 + \mu_3 - r)T - Z_T^2 - Z_T^3 + \sum_{j=2}^{3}(\beta_j + \frac{\sigma_j^2}{2})\gamma_T^j + \rho_{23}\sigma_2\sigma_3\sqrt{\gamma_T^2\gamma_T^3}}$$

$$\times \left(1 - N\left(\frac{K - \mu_1 T + Z_T^1 - \beta_I \gamma_T^1}{\sigma_1\sqrt{\gamma_T^1}}\right)\right)$$

$$= e^{(\mu_2 + \mu_3 - r)T - Z_T^2 - Z_T^3 + \sum_{j=2}^{3}(\beta_j + \frac{\sigma_j^2}{2})\gamma_T^j + \rho_{23}\sigma_2\sigma_3\sqrt{\gamma_T^2\gamma_T^3}}$$

$$\times N\left(\frac{\mu_1 T + \left(\beta_1 + \rho_{12}\sigma_1\sigma_2\sqrt{\frac{\gamma_T^2}{\gamma_T^1}} + \rho_{13}\sigma_1\sigma_3\sqrt{\frac{\gamma_T^3}{\gamma_T^1}}\right)\gamma_T^1 - K - Z_T^1}{\sigma_1\sqrt{\gamma_T^1}}\right). \quad (27)$$

Because $\rho_{12} = \rho_{13} = \rho_{23} = 0$, we get that

$$E\big(e^{-rT} S_T^3 S_T^2 I_{\{S_T^1 \geq K\}} | \gamma_T^1, Z_T^1, \gamma_T^2, Z_T^2, \gamma_T^3, Z_T^3\big)$$

$$= e^{(\mu_2 + \mu_3 - r)T - Z_T^2 - Z_T^3 + \sum_{j=2}^{3}(\beta_j + \frac{\sigma_j^2}{2})\gamma_T^j} N\left(\frac{\mu_1 T + \beta_1 \gamma_T^1 - K - Z_T^1}{\sigma_1\sqrt{\gamma_T^1}}\right)$$

$$= e^{(\mu_2 + \mu_3 - r)T - Z_T^2 - Z_T^3 + \left[\sum_{j=2}^{3} \kappa_j(\beta_j + \frac{\sigma_j^2}{2})\right]\gamma_T + \sum_{j=2}^{3}(\beta_j + \frac{\sigma_j^2}{2})\tilde{\kappa}_j \tilde{\gamma}_T^j}$$

$$\times e^{\left[\sum_{j=2}^{3} \kappa_{j1}(\beta_j + \frac{\sigma_j^2}{2})\right]\gamma_T^1} N\left(\frac{\mu_1 T + \beta_1 \gamma_T^1 - K - Z_T^1}{\sigma_1\sqrt{\gamma_T^1}}\right), \quad (28)$$

where $\kappa_j$, $\tilde{\kappa}_j$, $\kappa_{j1}$, $j = 2, 3$, are defined in (9).

Next, we pass to the computing of the conditional expectation

$$E\big(e^{-rT} S_T^3 S_T^2 I_{\{S_T^1 \geq K\}} | Z_T^1, Z_T^2, Z_T^3\big).$$

It is clear from (28) that we need to calculate the integral

$$I = \int_0^\infty x^\alpha e^{-(a_1 - b)x} N\left(h\sqrt{x} + \frac{p}{\sqrt{x}}\right) dx, \quad (29)$$

where $a_1$ is the parameter of $\gamma_t^1$ (see (8)),

$$\alpha = a_1 T - 1, \qquad p = \frac{\mu_1 T - K - Z_T^1}{\sigma_1}, \qquad b = \sum_{j=2}^{3} \kappa_{j1}\left(\beta_j + \frac{\sigma_j^2}{2}\right), \qquad h = \frac{\beta_1}{\sigma_1}.$$

Then

$$E\big(e^{-rT} S_T^3 S_T^2 I_{\{S_T^1 \geq K\}} | Z_T^1, Z_T^2, Z_T^3\big)$$

$$= e^{(\mu_2 + \mu_3 - r)T - Z_T^2 - Z_T^3} E\big(e^{\sum_{j=2}^{3}(\beta_j + \frac{\sigma_j^2}{2})\kappa_j \gamma_T}\big)$$

$$\times E\big(e^{(\beta_2 + \frac{\sigma_2^2}{2})\tilde{\kappa}_2 \tilde{\gamma}_T^2}\big) E\big(e^{(\beta_3 + \frac{\sigma_3^2}{2})\tilde{\kappa}_3 \tilde{\gamma}_T^3}\big) \frac{a_1^{a_1 T}}{\Gamma(a_1 T)} I. \quad (30)$$

Let us notice that the condition (5) is

$$E e^{\sum_{j=2}^{3} \beta_j \gamma_T^j + \frac{\sum_{j=2}^{3} \sigma_j^2 \gamma_T^j}{2}} < \infty$$

now since $\rho_{23} = 0$. Hence we have that

$$\mathrm{E}\big(e^{\sum_{j=2}^{3}(\beta_j+\frac{\sigma_j^2}{2})\kappa_j\gamma_T}\big)\mathrm{E}\big(e^{\sum_{j=2}^{3}(\beta_j+\frac{\sigma_j^2}{2})\kappa_{j1}\gamma_T^1}\big)$$
$$\times \mathrm{E}\big(e^{(\beta_2+\frac{\sigma_2^2}{2})\tilde{\kappa}_2\tilde{\gamma}_T^2}\big)\mathrm{E}\big(e^{(\beta_3+\frac{\sigma_3^2}{2})\tilde{\kappa}_3\tilde{\gamma}_T^3}\big) < \infty \quad (31)$$

and therefore $b < a_1$. And since $b < a_1$, we could apply to the integral (29) Cases 1–3 on pp. 207–212 of Ivanov and Ano [20]. If $p = 0$, then the identity

$$I = \frac{\Gamma(\alpha+\frac{3}{2})}{(a_1-b)^{\alpha+1}\sqrt{2\pi}}\bigg[\frac{\mathrm{B}(\frac{1}{2},\alpha+1)}{\sqrt{2}}$$
$$+ \frac{h}{\sqrt{a_1-b}}\mathrm{G}\bigg(\alpha+\frac{3}{2},\frac{1}{2},\frac{3}{2};-\frac{h^2}{2(a_1-b)}\bigg)\bigg] \quad (32)$$

is satisfied for $I$ defined in (29). When $p \neq 0$, we have that

$$I = \frac{|s|^{\alpha+\frac{1}{2}}e^s(1+q)^{\alpha+1}}{(a_1-b)^{\alpha+1}\sqrt{2\pi}}\bigg[\mathrm{B}(\alpha+1,1)\big(|s|\mathrm{M}_{\alpha+\frac{3}{2}}(|s|)$$
$$+ s\mathrm{M}_{\alpha+\frac{1}{2}}(|s|)\big)\mathrm{A}\bigg(\alpha+1,-\alpha,\alpha+2;\frac{1+q}{2},-s(1+q)\bigg)$$
$$- (1+q)s\mathrm{B}(\alpha+2,1)\mathrm{M}_{\alpha+\frac{1}{2}}(|s|)\mathrm{A}\bigg(\alpha+2,-\alpha,\alpha+3;\frac{1+q}{2},-s(1+q)\bigg)\bigg], \quad (33)$$

where

$$s = p\sqrt{h^2+2(a_1-b)} \quad \text{and} \quad q = \frac{h}{\sqrt{h^2+2(a_1-b)}}.$$

Set

$$\mathrm{DC}(n_1,n_2,n_3) = e^{-rT}\mathrm{E}\big(S_T^3 S_T^2 I_{\{S_T^1 \geq K\}}|N_T^1 = n_1, N_T^2 = n_2, N_T^3 = n_3\big).$$

Then we have that

$$\mathrm{DC}(n_1,n_2,n_3) = \mathrm{E}e^{X_T^2+X_T^3-\sum_{j=1}^{n_2}\xi_j^2-\sum_{j=1}^{n_3}\xi_j^3}I_{\{e^{X_T^1-\sum_{j=1}^{n_1}\xi_j^1}\geq K\}}$$
$$\geq \mathrm{E}e^{X_T^2+X_T^3-\sum_{j=1}^{\tilde{n}_2}\xi_j^2-\sum_{j=1}^{\tilde{n}_3}\xi_j^3}I_{\{e^{X_T^1-\sum_{j=1}^{\tilde{n}_1}\xi_j^1}\geq K\}} = \mathrm{DC}(\tilde{n}_1,\tilde{n}_2,\tilde{n}_3)$$

when $n_j \leq \tilde{n}_j, j = 1,2,3$. Therefore,

$$\sum_{n_1=0}^{N_1}\sum_{n_2=0}^{N_2}\sum_{n_3=0}^{N_3}\frac{\lambda_1^{n_1}\lambda_2^{n_2}\lambda_3^{n_3}T^{n_1+n_2+n_3}e^{-(\lambda_1+\lambda_2+\lambda_3)T}\mathrm{DC}(n_1,n_2,n_3)}{n_1!n_2!n_3!}$$
$$\leq \mathbb{DC}$$
$$\leq \sum_{n_1=0}^{N_1}\sum_{n_2=0}^{N_2}\sum_{n_3=0}^{N_3}\frac{\lambda_1^{n_1}\lambda_2^{n_2}\lambda_3^{n_3}T^{n_1+n_2+n_3}e^{-(\lambda_1+\lambda_2+\lambda_3)T}\mathrm{DC}(n_1,n_2,n_3)}{n_1!n_2!n_3!}$$

$$+ \text{DC}(n_1, n_2, n_3) \sum_{n_1=N_1+1}^{\infty} \sum_{n_2=N_2+1}^{\infty} \sum_{n_3=N_3+1}^{\infty} \frac{\lambda_1^{n_1}\lambda_2^{n_2}\lambda_3^{n_3}T^{n_1+n_2+n_3}}{n_1!n_2!n_3!e^{(\lambda_1+\lambda_2+\lambda_3)T}}$$

$$= \sum_{n_1=0}^{N_1}\sum_{n_2=0}^{N_2}\sum_{n_3=0}^{N_3} \frac{\lambda_1^{n_1}\lambda_2^{n_2}\lambda_3^{n_3}T^{n_1+n_2+n_3}e^{-(\lambda_1+\lambda_2+\lambda_3)T}\text{DC}(n_1,n_2,n_3)}{n_1!n_2!n_3!}$$

$$+ \text{DC}(n_1,n_2,n_3)\left(1 - \sum_{n_1=0}^{N_1}\frac{\lambda_1^{n_1}e^{-\lambda_1 T}}{n_1!}\right)$$

$$\times \left(1 - \sum_{n_2=0}^{N_2}\frac{\lambda_2^{n_2}e^{-\lambda_2 T}}{n_2!}\right)\left(1 - \sum_{n_3=0}^{N_3}\frac{\lambda_3^{n_3}e^{-\lambda_3 T}}{n_3!}\right). \qquad (34)$$

The result of Theorem 1 follows from (34), where the functions $\text{DC}(n_1,n_2,n_3)$ are computed with respect to (30) using (7) and (32)–(33). □

**Proof of Theorem 2.** Since $\rho_{12} = \rho_{13} = 0$ and $\gamma_T^3 \equiv \gamma_T^2$, we get using (27) that

$$\text{E}\big(e^{-rT}S_T^3 S_T^2 I_{\{S_T^1 \geq K\}} | \gamma_T^1, Z_T^1, \gamma_T^2, Z_T^2, \gamma_T^3, Z_T^3\big)$$
$$= e^{(\mu_2+\mu_3-r)T - Z_T^2 - Z_T^3 + \kappa_2\left[\sum_{j=2}^{3}(\beta_j + \frac{\sigma_j^2}{2}) + \rho_{23}\sigma_2\sigma_3\right]\gamma_T}$$
$$\times e^{\kappa_{21}\left[\sum_{j=2}^{3}(\beta_j + \frac{\sigma_j^2}{2}) + \rho_{23}\sigma_2\sigma_3\right]\gamma_T^1} \text{N}\left(\frac{\mu_1 T + \beta_1 \gamma_T^1 - K - Z_T^1}{\sigma_1\sqrt{\gamma_T^1}}\right). \qquad (35)$$

To get

$$\text{E}\big(e^{-rT}S_T^3 S_T^2 I_{\{S_T^1 \geq K\}} | Z_T^1, Z_T^2, Z_T^3\big),$$

we need to calculate the integral $I$ (29) with

$$\alpha = a_1 T - 1, \qquad p = \frac{\mu_1 T - K - Z_T^1}{\sigma_1},$$

$$b = \kappa_{21}\left[\sum_{j=2}^{3}\left(\beta_j + \frac{\sigma_j^2}{2}\right) + \rho_{23}\sigma_2\sigma_3\right], \qquad h = \frac{\beta_1}{\sigma_1}.$$

Then

$$\text{E}\big(e^{-rT}S_T^3 S_T^2 I_{\{S_T^1 \geq K\}} | Z_T^1, Z_T^2, Z_T^3\big) = e^{(\mu_2+\mu_3-r)T - Z_T^2 - Z_T^3}$$
$$\times \text{E}\big(e^{\kappa_2\left[\sum_{j=2}^{3}(\beta_j+\frac{\sigma_j^2}{2}) + \rho_{23}\sigma_2\sigma_3\right]\gamma_T}\big)\frac{a_1^{a_1 T}}{\Gamma(a_1 T)}I, \qquad (36)$$

where $I$ is calculated by (32)–(33).

Under the conditions of Theorem 2, (3) has the form

$$\text{E}e^{\gamma_T^2 \sum_{j=2}^{3}\beta_j + \frac{\sum_{j=2}^{3}\sigma_j^2 + 2\rho_{23}\sigma_2\sigma_3}{2}\gamma_T^2} < \infty. \qquad (37)$$

Therefore,

$$\mathrm{E}\bigl(e^{\kappa_2\bigl[\sum_{j=2}^{3}(\beta_j+\frac{\sigma_j^2}{2})+\rho_{23}\sigma_2\sigma_3\bigr]\gamma_T}\bigr)\mathrm{E}\bigl(e^{\kappa_{21}\bigl[\sum_{j=2}^{3}(\beta_j+\frac{\sigma_j^2}{2})+\rho_{23}\sigma_2\sigma_3\bigr]\gamma_T^1}\bigr) < \infty$$

and $b < a_1$. The result of Theorem 2 comes from (36) analogously to the result of Theorem 1. □

**Proof of Theorem 3.** Because $\gamma_T^3 = \gamma_T^2 = \gamma_T^1$, we have from (27) that

$$\mathrm{E}\bigl(e^{-rT}S_T^3 S_T^2 I_{\{S_T^1 \geq K\}}|\gamma_T^1, Z_T^1, \gamma_T^2, Z_T^2, \gamma_T^3, Z_T^3\bigr)$$
$$= e^{(\mu_2+\mu_3-r)T - Z_T^2 - Z_T^3 + \bigl[\sum_{j=2}^{3}(\beta_j+\frac{\sigma_j^2}{2})+\rho_{23}\sigma_2\sigma_3\bigr]\gamma_T^1}$$
$$\times \mathrm{N}\left(\frac{\mu_1 T + (\beta_1 + \rho_{12}\sigma_1\sigma_2 + \rho_{13}\sigma_1\sigma_3)\gamma_T^1 - K - Z_T^1}{\sigma_1\sqrt{\gamma_T^1}}\right). \quad (38)$$

Hence

$$\mathrm{E}\bigl(e^{-rT}S_T^3 S_T^2 I_{\{S_T^1 \geq K\}}|Z_T^1, Z_T^2, Z_T^3\bigr) = e^{(\mu_2+\mu_3-r)T - Z_T^2 - Z_T^3}\frac{a_1^{a_1 T}}{\Gamma(a_1 T)}I, \quad (39)$$

where $I$ is defined in (29) and computed by (32)–(33) with the same $\alpha$ and $p$ as in the proof of Theorem 2,

$$b = \sum_{j=2}^{3}\left(\beta_j + \frac{\sigma_j^2}{2}\right) + \rho_{23}\sigma_2\sigma_3, \qquad h = \frac{\beta_1 + \rho_{12}\sigma_1\sigma_2 + \rho_{13}\sigma_1\sigma_3}{\sigma_1}.$$

The condition (3) in Theorem 3 has the form

$$\mathrm{E}\bigl(e^{\bigl[\sum_{j=2}^{3}(\beta_j+\frac{\sigma_j^2}{2})+\rho_{23}\sigma_2\sigma_3\bigr]\gamma_T^1}\bigr) < \infty$$

now and hence $b < a_1$. The result of Theorem 3 is derived from (39) using (7) and (32)–(34) from the proof of Theorem 1. □

**Proof of Theorem 4.** Keeping in mind the conditions of Theorem 4, one could observe from (27) that

$$\mathrm{E}\bigl(e^{-rT}S_T^3 S_T^2 I_{\{S_T^1 \geq K\}}|\gamma_T^1, Z_T^1, \gamma_T^2, Z_T^2, \gamma_T^3, Z_T^3\bigr)$$
$$= e^{(\mu_2+\mu_3-r)T - Z_T^2 - Z_T^3 + \tilde{\kappa}_2(\beta_2+\frac{\sigma_2^2}{2})\tilde{\gamma}_T^2 + \bigl[\beta_3+\frac{\sigma_3^2}{2}+\kappa_{21}(\beta_2+\frac{\sigma_2^2}{2})\bigr]\gamma_T^1}$$
$$\times \mathrm{N}\left(\frac{\mu_1 T + (\beta_1 + \sigma_1\sigma_3)\gamma_T^1 - K - Z_T^1}{\sigma_1\sqrt{\gamma_T^1}}\right). \quad (40)$$

Therefore

$$\mathrm{E}\bigl(e^{-rT}S_T^3 S_T^2 I_{\{S_T^1 \geq K\}}|Z_T^1, Z_T^2, Z_T^3\bigr)$$
$$= e^{(\mu_2+\mu_3-r)T - Z_T^2 - Z_T^3}\mathrm{E}\bigl(e^{\tilde{\kappa}_2(\beta_2+\frac{\sigma_2^2}{2})\tilde{\gamma}_T^2}\bigr)\frac{a_1^{a_1 T}}{\Gamma(a_1 T)}I, \quad (41)$$

where $I$ is defined in (29) and computed by (32)–(33) with the same $\alpha$ and $p$ as in the proof of Theorem 2,

$$b = \beta_3 + \frac{\sigma_3^2}{2} + \kappa_{21}\left(\beta_2 + \frac{\sigma_2^2}{2}\right), \qquad h = \frac{\beta_1 + \sigma_1\sigma_3}{\sigma_1}.$$

The condition (3) has here the form

$$\mathrm{E}\big(e^{(\beta_3 + \frac{\sigma_3^2}{2})\gamma_T^1 + (\beta_2 + \frac{\sigma_2^2}{2})\gamma_T^2}\big) < \infty. \tag{42}$$

Therefore,

$$\mathrm{E}\big(e^{\left[\beta_3 + \frac{\sigma_3^2}{2} + \kappa_{21}(\beta_2 + \frac{\sigma_2^2}{2})\right]\gamma_T^1}\big)\mathrm{E}\big(e^{\tilde{\kappa}_2(\beta_2 + \frac{\sigma_2^2}{2})\tilde{\gamma}_T^2}\big) < \infty \tag{43}$$

and hence $b < a_1$. It means that we can exploit here the results of Ivanov and Ano [20] and obtain the result of Theorem 4 from (41) in the same way as it is made in the proof of Theorem 1 in (32)–(34). □

**Proof of Theorem 5.** We have similarly to (28) that

$$\mathrm{E}\big(e^{-rT}S_T^3 S_T^2 I_{\{S_T^1 \geq K\}}\big|\varkappa_T^1, Z_T^1, \varkappa_T^2, Z_T^2, \varkappa_T^3, Z_T^3\big)$$

$$= e^{(\mu_2 + \mu_3 - r)T - Z_T^2 - Z_T^3 + \left[\sum_{j=2}^3 \kappa_j(\beta_j + \frac{\sigma_j^2}{2})\right]\varkappa_T + \sum_{j=2}^3 (\beta_j + \frac{\sigma_j^2}{2})\tilde{\kappa}_j \tilde{\varkappa}_T^j}$$

$$\times e^{\left[\sum_{j=2}^3 \kappa_{j1}(\beta_j + \frac{\sigma_j^2}{2})\right]\varkappa_T^1} \mathrm{N}\left(\frac{\mu_1 T + \beta_1 \varkappa_T^1 - K - Z_T^1}{\sigma_1 \sqrt{\varkappa_T^1}}\right).$$

Since

$$\frac{\phi_1^2}{2} = \sum_{j=2}^3 \kappa_{j1}\left(\beta_j + \frac{\sigma_j^2}{2}\right)$$

with respect to the conditions of Theorem 5, one can notice that

$$\mathrm{E}\big(e^{-rT}S_T^3 S_T^2 I_{\{S_T^1 \geq K\}}|Z_T^1, Z_T^2, Z_T^3\big) = e^{(\mu_2 + \mu_3 - r)T - Z_T^2 - Z_T^3}$$

$$\times \mathrm{E}\big(e^{(\beta_2 + \frac{\sigma_2^2}{2})\tilde{\kappa}_2 \tilde{\varkappa}_T^2}\big)\mathrm{E}\big(e^{(\beta_3 + \frac{\sigma_3^2}{2})\tilde{\kappa}_3 \tilde{\varkappa}_T^3}\big)$$

$$\times \mathrm{E}\big(e^{\left[\sum_{j=2}^3 \kappa_j(\beta_j + \frac{\sigma_j^2}{2})\right]\varkappa_T}\big)\frac{\phi_1 T e^{\phi_1^2 T}}{\sqrt{2\pi}} J,$$

where

$$J = \int_0^\infty x^{-\frac{3}{2}} e^{-\frac{(\phi_1 T)^2}{2x}} \mathrm{N}\left(\frac{\mu_1 T + \beta_1 x - K - Z_T^1}{\sigma_1 \sqrt{x}}\right) dx$$

$$= \frac{\sqrt{2}}{\phi_1 T}\int_0^\infty x^{-\frac{3}{2}} e^{-\frac{1}{x}} \mathrm{N}\left(h\sqrt{x} + \frac{p}{\sqrt{x}}\right) \tag{44}$$

with

$$h = \frac{\beta_1 \phi_1 T}{\sigma_1 \sqrt{2}} \quad \text{and} \quad p = \frac{(\mu_1 T - K - Z_T^1)\sqrt{2}}{\sigma_1 \phi_1 T}.$$

If $\beta_1 \neq 0$, it is easy to see that the integral (44) is quite the same as the integral (4.1) in Ivanov and Temnov [21]. Therefore, we get from (4.3)–(4.6) of Ivanov and Temnov [21] that if $\beta_1 \neq 0$ then

$$J = \frac{1}{\phi_1 T \sqrt{2}}(J_1 + J_2), \tag{45}$$

where

$$J_1 = |\varsigma|(q+1)^{-\frac{1}{2}} \exp(|\varsigma|) M_1(|\varsigma|) \Upsilon_1$$

and

$$J_2 = |\varsigma|(q+1)^{-\frac{1}{2}} \exp(|\varsigma|) M_0(|\varsigma|) (\Upsilon_1 - (q+1)\Upsilon_2)$$

with

$$\Upsilon_1 = B\left(\frac{1}{2}, 1\right) A\left(\frac{1}{2}, \frac{1}{2}, \frac{3}{2}; \frac{q+1}{2}, -|\varsigma|(q+1)\right)$$

and

$$\Upsilon_2 = B\left(\frac{3}{2}, 1\right) A\left(\frac{3}{2}, \frac{1}{2}, \frac{5}{2}; \frac{q+1}{2}, -|\varsigma|(q+1)\right),$$

where

$$\varsigma = h\sqrt{p^2 + 2} \quad \text{and} \quad q = \frac{p}{\sqrt{p^2 + 2}}.$$

When $\beta_1 = 0$, we have that

$$\begin{aligned}
J &= \frac{\sqrt{2}}{\phi_1 T} \int_0^\infty x^{-\frac{3}{2}} e^{-\frac{1}{x}} N\left(\frac{p}{\sqrt{x}}\right) dx \\
&= \frac{\sqrt{2}}{\phi_1 T} \int_0^\infty x^{-\frac{3}{2}} e^{-\frac{1}{x}} \left(\int_{-\infty}^{\frac{p}{\sqrt{x}}} \frac{1}{\sqrt{2\pi}} e^{-\frac{y^2}{2}} dy\right) dx \\
&= \frac{\sqrt{2}}{\phi_1 T} \int_0^\infty x^{-\frac{3}{2}} e^{-\frac{1}{x}} \left(\int_{-\infty}^p \frac{1}{\sqrt{2\pi x}} e^{-\frac{y^2}{2x}} dy\right) dx \\
&= \frac{1}{\phi_1 T \sqrt{\pi}} \int_{-\infty}^p \left(\int_0^\infty x^{-2} e^{-\frac{1}{x} - \frac{y^2}{2x}} dx\right) dy. \tag{46}
\end{aligned}$$

Let us notice that the Fubini theorem can be applied to $J$ since the double integral

$$\int_0^\infty \int_{-\infty}^p x^{-2} e^{-\frac{1}{x} - \frac{y^2}{2x}} dy dx$$

is an integral of constant sign function and because the Fubini theorem is applicable to

$$\int_0^n \int_{-n}^p x^{-2} e^{-\frac{1}{x} - \frac{y^2}{2x}} dy dx$$

for any $n \in \mathbb{N}$ as the integrand is continuous. Because

$$\int_0^\infty x^{-2} e^{-\frac{1}{x} - \frac{y^2}{2x}} dx = \left(1 + \frac{y^2}{2}\right)^{-1} \int_0^\infty \left(1 + \frac{y^2}{2}\right) x^{-2} e^{-\frac{1}{x} - \frac{y^2}{2x}} dx$$

$$= \left(1 + \frac{y^2}{2}\right)^{-1} \int_0^\infty de^{-(1+\frac{y^2}{2})\frac{1}{x}} = \left(1 + \frac{y^2}{2}\right)^{-1},$$

it follows from (46) that

$$J = \frac{1}{\phi_1 T \sqrt{\pi}} \int_{-\infty}^p \left(1 + \frac{y^2}{2}\right)^{-1} dy = \frac{\sqrt{2}}{\phi_1 T \sqrt{\pi}} \int_{-\infty}^{\frac{p}{\sqrt{2}}} (1 + y^2)^{-1} dy$$

$$= \frac{\sqrt{2}}{\phi_1 T \sqrt{\pi}} \left(\frac{\pi}{2} + \mathrm{sign}\, p \arctan \frac{|p|}{\sqrt{2}}\right) \tag{47}$$

if $\beta_1 = 0$.

If $a^2 > 2A$, it follows from (19) that

$$\mathrm{E} e^{A \varkappa_t} = \frac{\phi t e^{a \phi t}}{\sqrt{2\pi}} \int_0^\infty x^{-\frac{3}{2}} e^{-\frac{1}{2}\left((a^2 - 2A)x + \frac{(\phi t)^2}{x}\right)} dx = e^{\phi t (a - \sqrt{a^2 - 2A})}. \tag{48}$$

When $a^2 = 2A$, the expectation

$$\mathrm{E} e^{A \varkappa_t} = \frac{\phi t e^{a \phi t}}{\sqrt{2\pi}} \int_0^\infty x^{-\frac{3}{2}} e^{-\frac{(\phi t)^2}{2x}} dx = \frac{\phi t e^{a \phi t}}{\sqrt{2\pi}} \int_{-\infty}^0 |x|^{-\frac{1}{2}} e^{\frac{(\phi t)^2 x}{2}} dx$$

$$= \frac{\phi t e^{a \phi t}}{\sqrt{2\pi}} \int_0^\infty x^{-\frac{1}{2}} e^{-\frac{(\phi t)^2 x}{2}} dx = \frac{e^{a \phi t} \Gamma(\frac{1}{2})}{\sqrt{\pi}} = e^{a \phi t}. \tag{49}$$

The condition (3) has the form (31) here. Therefore,

$$\phi^2 \geq 2 \sum_{j=2}^3 \kappa_j \left(\beta_j + \frac{\sigma_j^2}{2}\right) \quad \text{and} \quad \tilde{\phi}_j^2 \geq 2\left(\beta_j + \frac{\sigma_j^2}{2}\right) \tilde{\kappa}_j, \quad j = 2, 3.$$

Hence the result of Theorem 5 comes from (45), (47) and (48)–(49). □

**Proof of Theorem 6.** Since

$$\phi_1^2 = 2\kappa_{21} \left(\sum_{j=2}^3 \left(\beta_j + \frac{\sigma_j^2}{2}\right) + \rho_{23} \sigma_2 \sigma_3\right),$$

we get using (35) that

$$\mathrm{E}\left(e^{-rT} S_T^3 S_T^2 I_{\{S_T^1 \geq K\}} | Z_T^1, Z_T^2, Z_T^3\right)$$

$$= e^{(\mu_2 + \mu_3 - r)T - Z_T^2 - Z_T^3} \mathrm{E}\left(e^{\kappa_2 \left[\sum_{j=2}^3 (\beta_j + \frac{\sigma_j^2}{2}) + \rho_{23} \sigma_2 \sigma_3\right] \varkappa_T}\right) \frac{\phi_1 T e^{\phi_1^2 T}}{\sqrt{2\pi}} J,$$

where $J$ is defined (44). The condition (3) has the form (37) here. Therefore,

$$\phi^2 \geq 2\kappa_2 \left[\sum_{j=2}^3 \left(\beta_j + \frac{\sigma_j^2}{2}\right) + \rho_{23} \sigma_2 \sigma_3\right]$$

and hence

$$\mathrm{E}\bigl(e^{\kappa_2\bigl[\sum_{j=2}^{3}(\beta_j+\frac{\sigma_j^2}{2})+\rho_{23}\sigma_2\sigma_3\bigr]\varkappa_T}\bigr) < \infty$$

and can be computed by (48) and (49). □

**Proof of Theorem 7.** We have from (38) and the condition

$$\phi_1^2 = 2\sum_{j=2}^{3}\left(\beta_j + \frac{\sigma_j^2}{2}\right) + \rho_{23}\sigma_2\sigma_3$$

that

$$\mathrm{E}\bigl(e^{-rT}S_T^3 S_T^2 I_{\{S_T^1 \geq K\}} | Z_T^1, Z_T^2, Z_T^3\bigr) = e^{(\mu_2+\mu_3-r)T - Z_T^2 - Z_T^3} \frac{\phi_1 T e^{\phi_1^2 T}}{\sqrt{2\pi}} J$$

with

$$J = \frac{\sqrt{2}}{\phi_1 T} \int_0^\infty x^{-\frac{3}{2}} e^{-\frac{1}{x}} \mathrm{N}\left(h\sqrt{x} + \frac{p}{\sqrt{x}}\right), \tag{50}$$

where

$$h = \frac{(\beta_1 + \rho_{12}\sigma_1\sigma_2 + \rho_{13}\sigma_1\sigma_3)\phi_1 T}{\sigma_1\sqrt{2}} \quad \text{and} \quad p = \frac{(\mu_1 T - K - Z_T^1)\sqrt{2}}{\sigma_1 \phi_1 T}.$$

Hence $J$ is determined by (45) if $\beta_1 + \rho_{12}\sigma_1\sigma_2 + \rho_{13}\sigma_1\sigma_3 \neq 0$ and by (47) when $\beta_1 + \rho_{12}\sigma_1\sigma_2 + \rho_{13}\sigma_1\sigma_3 = 0$. □

**Proof of Theorem 8.** The condition (3) has the form (42)–(43) here. Therefore, $\tilde{\phi}_2^2 \geq \tilde{\kappa}_2(2\beta_2 + \sigma_2^2)$ and we get from (40) and the condition $\phi_1^2 = 2\beta_3 + \sigma_3^2 + \kappa_{21}(2\beta_2 + \sigma_2^2)$ that

$$\mathrm{E}\bigl(e^{-rT}S_T^3 S_T^2 I_{\{S_T^1 \geq K\}} | Z_T^1, Z_T^2, Z_T^3\bigr)$$
$$= e^{(\mu_2+\mu_3-r)T - Z_T^2 - Z_T^3} \mathrm{E}\bigl(e^{\tilde{\kappa}_2(\beta_2+\frac{\sigma_2^2}{2})\tilde{\varkappa}_T^2}\bigr) \frac{\phi_1 T e^{\phi_1^2 T}}{\sqrt{2\pi}} J,$$

where $J$ is defined in (50) with

$$h = \frac{(\beta_1 + \sigma_1\sigma_3)\phi_1 T}{\sigma_1\sqrt{2}} \quad \text{and} \quad p = \frac{(\mu_1 T - K - Z_T^1)\sqrt{2}}{\sigma_1 \phi_1 T}.$$

Therefore, $J$ is computed by (45) if $\beta_1 + \sigma_1\sigma_3 \neq 0$ and by (47) when $\beta_1 + \sigma_1\sigma_3 = 0$. □

106    R.V. Ivanov, K. Ano[19] Ivanov, R.V.: Option pricing in the variance-gamma model under the drift jump. Int. J. Theor. Appl. Finance **21**(4), 1–19 (2018). MR3817272. 10.1142/S0219024918500188

[20] Ivanov, R.V., Ano, K.: On exact pricing of fx options in multivariate time-changed Lévy models. Rev. Deriv. Res. **19**(3), 201–216 (2016). MR3611537. 10.1007/s11009-015-9461-8

[21] Ivanov, R.V., Temnov, G.: Truncated moment-generating functions of the nig process and their applications. Stoch. Dyn. **17**(5), 1–12 (2017). MR3669415. 10.1142/S0219493717500393

[22] Kallsen, J., Shiryaev, A.N.: The cumulant process and Esscher's change of measure. Finance Stoch. **6**, 397–428 (2002). MR1932378. 10.1007/s007800200069

[23] Koponen, I.: Analytic approach to the problem of convergence of truncated Lévy flights towards the Gaussian stochastic process. Phys. Rev. E **52**, 1197–1199 (1995)

[24] Küchler, U., Tappe, S.: Bilateral gamma distributions and processes in financial mathematics. Stoch. Process. Appl. **118**(2), 261–283 (2008). MR2376902. 10.1016/j.spa.2007.04.006

[25] Küchler, U., Tappe, S.: Tempered stable distributions and processes. Stoch. Process. Appl. **123**(12), 4256–4293 (2013). MR3096354. 10.1016/j.spa.2013.06.012

[26] Linders, D., Stassen, B.: The multivariate variance gamma model: basket option pricing and calibration. Quant. Finance **16**(4), 555–572 (2016). MR3473973. 10.1080/14697688.2015.1043934

[27] Luciano, E., Schoutens, W.: A multivariate jump-driven financial asset model. Quant. Finance **6**(5), 385–402 (2016). MR2261218. 10.1080/14697680600806275

[28] Luciano, E., Semeraro, P.: Multivariate time changes for Lévy asset models: Characterization and calibration. J. Comput. Appl. Math. **233**, 1937–1953 (2010). MR2564029. 10.1016/j.cam.2009.08.119

[29] Luciano, E., Marena, M., Semeraro, P.: Dependence calibration and portfolio fit with factor-based subordinators. Quant. Finance **16**(7), 1–16 (2016). MR3516138. 10.1080/14697688.2015.1114661

[30] Madan, D., Milne, F.: Option pricing with vg martingale components. Math. Finance **1**(4), 39–55 (1991)

[31] Madan, D., Seneta, E.: The variance gamma (v.g.) model for share market returns. J. Bus. **63**, 511–524 (1991)

[32] Madan, D., Yor, M.: Representing the cgmy and Meixner Lévy processes as time changed Brownian motions. J. Comput. Finance **12**(1), 27–47 (2008). MR2504899. 10.21314/JCF.2008.181

[33] Madan, D., Carr, P., Chang, E.: The variance gamma process and option pricing. Eur. Finance Rev. **2**, 79–105 (1998)

[34] Moosbrucker, T.: Explaining the correlation smile using variance gamma distributions. J. Fixed Income **16**(1), 71–87 (2006)

[35] Mozumder, S., Sorwar, G., Dowd, K.: Revisiting variance gamma pricing: An application to S&P500 index options. Int. J. Financ. Eng. **2**(2), 1–24 (2015). MR3454654. 10.1142/S242478631550022X

[36] Rathgeber, A., Stadler, J., Stöckl, S.: Modeling share returns – an empirical study on the variance gamma model. J. Econ. Finance **40**(4), 653–682 (2016)